\documentclass[12pt]{amsart}

\usepackage{amsmath,amssymb,amsthm}

\usepackage{fancyhdr}
\usepackage{amsfonts,graphicx}

\theoremstyle{plain}
\newtheorem{Theo}{Theorem}[section]
\newtheorem{lem}[Theo]{Lemma}
\newtheorem{prop}[Theo]{Proposition}

\theoremstyle{plain}
\theoremstyle{definition}

\newtheorem*{nota}{Notation}
\theoremstyle{remark}

\newtheorem*{rema*}{Remarks}

\newcommand{\ZZ}{\mathbb{Z}}  
  
\newcommand{\NN}{\mathbb{N}}

\newcommand{\RR}{\mathbb{R}}

\newcommand{\diver}[1]{ \textnormal{div}\hspace{0.07cm} #1} 

\date{}

\begin{document}

\title[Global  existence for the  axisymmetric Euler equations]
{On the global existence for the axisymmetric Euler equations}
\author[H. Abidi]{Hammadi Abidi}
\address{IRMAR, Universit\'e de Rennes 1\\ Campus de
Beaulieu\\ 35~042 Rennes cedex\\ France}
\email{hamadi.abidi@univ-rennes1.fr}
\author[T. Hmidi]{Taoufik Hmidi}
\address{IRMAR, Universit\'e de Rennes 1\\ Campus de
Beaulieu\\ 35~042 Rennes cedex\\ France}
\email{thmidi@univ-rennes1.fr}
\author[S. Keraani]{Sahbi Keraani}
\address{IRMAR, Universit\'e de Rennes 1\\ Campus de
Beaulieu\\ 35~042 Rennes cedex\\ France}
\email{sahbi.keraani@univ-rennes1.fr}
\begin{abstract}
This paper deals with the global well-posedness of  the $3$D
axisymmetric Euler equations for  initial data lying in some  critical Besov spaces.
\end{abstract}
\maketitle

\section{Introduction}
The motions of homogeneous inviscid incompressible fluid flows in $\RR^3$ are governed by the
Euler system
$$
\mbox{(E)}\left\lbrace
\begin{array}{l}
\partial_t u+(u\cdot\nabla)u+\nabla \pi=0,\\
\diver u=0,\\
{u}_{| t=0}=u_0.
\end{array}
\right.
$$ 
Here,  $u=u(t,x)\in \RR^3$ denotes the   velocity of the fluid, $\pi=\pi(t,x)$  the  scalar pressure and $u_0$  is the initial velocity satisfying $\diver u_0=0$.

The local well-posedness theory of the system (E) seems to be in a satisfactory state and several  results are obtained by numerous authors  in many standard function spaces. In \cite{Kato}, Kato proved the local existence  and uniqueness for initial data $u_0\in H^s(\RR^3)$ with $s>5/2$ and \mbox{Chemin \cite{Ch1}} gave similar results for initial data lying in H\"olderian \mbox{spaces $C^r$} 
with $r>1.$ Other local results are recently obtained by \mbox{Chae \cite{Chae}} in critical Besov space 
$B_{p,1}^{\frac{3}{p}+1},$ with $p\in]1,\infty[$ and by Pak and Park \cite{Hee} for the space $B_{\infty,1}^1$. Nevertheless the question of global existence is still open
and continues to be  one of the leading problem in mathematical fluid mechanics.  On the other hand, there are many criteria of finite time blowup. One of them is  the  BKM criterion  \cite{Beale} which ensures that the developement of finite time singularities is related to the divergence of the $L^\infty$ norm of the vorticity near the maximal existence time. A direct consequence of this result is the global existence in time  of  smooth two-dimensional Euler solutions since the vorticity is only convected and then not grow. We emphasize that new geometric blowup criteria are recently discovered by Constantin, Fefferman and Majda \cite{CFM}. 

We recall that the vorticity  $\omega=\nabla\times u$  satisfies the equation 
$$
\partial_t \omega+(u\cdot\nabla)\omega-(\omega \cdot\nabla)u=0.
$$
 The main 
difficulty for establishing global regularity is to understand how the vortex stretching term
$(\omega \cdot\nabla)u$ affects the dynamic of the fluid. While global existence is not proved for arbitrary initial smooth data, there are partial results in the case of the so-called axisymmetric flows without swirl.
 By an axisymmetric solution without swirl of 
the Euler  system we mean a solution of the form
$$
u(x, t) = u^r(r, z , t)e_r + u^z (r, z , t)e_z,\quad x=(x_{1},x_{2},z),\quad r=\sqrt{x_{1}^2+x_{2}^2},
$$ 
where  $\big(e_r, e_{\theta} , e_z\big)$ is the cylindrical basis of $\RR^3$. 
The corresponding  vorticity has the form
$$
\omega=(\partial_zu^r-\partial_ru^z) e_{\theta}
$$
and satisfies 
$$
\partial_t \omega +(u^r\partial_r+u^z\partial_z)\omega =\frac{u^r}{r}\omega.
$$
Since $u^{\theta}=0$, then the equation can be rewritten as
  \begin{equation}
 \label{tourbillon}
\partial_t \omega +(u\cdot\nabla)\omega =\frac{u^r}{r}\omega.
\end{equation}
One of the main property of the axisymmetric flows is the preservation of the  quantity $\alpha:=\omega/r$ along the flow, that is 
\begin{equation}
\label{equation_importante}
\partial_t\alpha+(u\cdot\nabla)\alpha=0.
\end{equation}
This fact was crucial to  Ukhovskii and Iudovich \cite{Ukhovskii}
to prove the global existence for initial data in  $H^s$ with $s>{7\over 2}$.
This result was improved by  Shirota and Yanagisawa \cite{Taira} where the condition $s>{7\over 2}$ is relaxed \mbox{to  $s>{5\over 2}$.}
Their proof is based on the boundedness of the quantity $\frac{u^r}{r}$ by using Biot-Savart law. We mention also the reference  \cite{Saint} where  similar results  are given in different function spaces.

Our main goal in this paper is to prove the global existence and uniqueness for more rough initial data lying in  some critical Besov spaces. 
Before going further into details, we fix some notations.
\begin{nota}
i) For $p_{1},p_{2}\in[1,+\infty]$ and $s\in\RR$ we denote respectively by $B^s_{p_1,p_2}$ and $L^{p_{1},p_{2}}$ the Besov  and   Lorentz spaces (see  next section for more details about these spaces).

ii) The space  $ \widetilde B^{0}_{\infty,1}$ is defined by the set of $v\in \mathcal{S}'(\RR^3)$ such that
$$
\Vert v \Vert_{ \widetilde B^{0}_{\infty,1}}:= \sum_{q\geq -1}
\| v-S_{q}v\|_{L^\infty}<\infty.
$$
See also next section for the definition of the operator $S_q$.\\
iii) For every $k\in \Bbb N$,  $\Phi_k$ denotes a real valued function having the form   
$$
\Phi_k(t)=  C_{0}\underbrace{ \exp(...\exp  }_{k\,times}(C_0t)...),
$$  
where $C_0$ is a constant which depends polynomially on the   norm of the initial data $u_0$.
\end{nota}
The main result of this paper can  be stated as follows.
\begin{Theo}\label{thm0}  Let $u_{0}$ be an axisymmetric divergence free vector field belonging to  $ B_{p,1}^{\frac{3}{p}+1},$ with $p\in]1,\infty].$  We assume in addition that its vorticity satisfies $\omega_{0}\in\widetilde B_{\infty,1}^0$ and $\frac{\omega_{0}}{r}\in L^{3,1}$. Then there exists a unique global solution
$
u\in {\mathcal C}(\RR_+;\,  B^{1+{3\over p}}_{p,1})$ 
to  the system ${\rm (E)}$. 

Besides, there exists a function $\Phi_{4}$ such that
  $$
\forall t\in\RR_{+},\,\quad\|u(t)\|_{ B^{1+{3\over p}}_{p,1}}
\leq \Phi_4(t).
$$

\end{Theo}
\begin{rema*}\label{RK1}

 i) The space $\widetilde B_{\infty,1}^0$ is a subspace of $B_{\infty,1}^0$ and contains the space
$$
\Big\{u\in \mathcal{S}' , \sum_{q\geq-1}(q+2)\|\Delta_{q}v\|_{L^\infty}<\infty\Big\}\cdot
$$
ii) We mention that for $p<3$ the condition $\frac{\omega_{0}}{r}\in L^{3,1}$ is automatically derived from $u_{0}\in B_{p,1}^{\frac{3}{p}+1},$ see Proposition \ref{lorentz} for more further precisions.
\end{rema*}
The proof uses two crucial estimates, the first  one is the boundedness of the vorticity for every time which is obtained from Biot-Savart law and the use of Lorentz spaces, see Proposition \ref{ur0}. Unfortunately this is not sufficient to show global existence because we do not know whether the BKM criterion  works in the critical spaces or not. Thus we are led to establish a second new estimate for the vorticity in  Besov \mbox{space $B_{\infty,1}^0$} (see Proposition \ref{ur1}). This allows us to bound for every time the Lipschitz norm of the velocity which is  sufficient to prove global existence.

The rest of this paper is organized as follows. In section $ 2$   we recall some function spaces and gather some  preliminary estimates. The  proof  of  
\mbox{Theorem \ref{thm0}} is  given in section 3.
 
\section{Notations and preliminaries} 
Throughout this paper, $C$ stands for some real positive constant which may be different in each occurrence. We shall sometimes alternatively use the notation $X\lesssim Y$ for an inequality of type $X\leq CY$.
We shall also use the  notation
$$
U(t):=\int_0^t\|\nabla u(\tau)\|_{L^\infty}d\tau.
$$

Let us start with  a classical dyadic decomposition of the full space (see for instance \cite{Ch1}):
there exist two radially   functions  $\chi\in \mathcal{D}(\RR^3)$ and  
$\varphi\in\mathcal{D}(\RR^3\backslash{\{0\}})$ such that
\begin{itemize}
\item[\textnormal{i)}]
$\displaystyle{\chi(\xi)+\sum_{q\geq0}\varphi(2^{-q}\xi)=1},$
\item[\textnormal{ii)}]
$ \textnormal{supp }\varphi(2^{-p}\cdot)\cap
\textnormal{supp }\varphi(2^{-q}\cdot)=\varnothing,$ if  $|p-q|\geq 2$,\\

\item[\textnormal{iii)}]
$\displaystyle{q\geq1\Rightarrow \textnormal{supp}\chi\cap \textnormal{supp }\varphi(2^{-q})=\varnothing}$.
\end{itemize}
For every $u\in{\mathcal S}'$ one defines the nonhomogeneous Littlewood-Paley operators
$$
\Delta_{-1}u=\chi(\hbox{D})u;\, \forall
q\in\NN,\;\Delta_qu=\varphi(2^{-q}\hbox{D})u\; \quad\hbox{and}\;
S_qu=\sum_{-1\leq j\leq q-1}\Delta_{j}u.
$$
The homogeneous operators are defined as follows
$$
\forall q\in\ZZ,\,\quad\dot{\Delta}_{q}u=\varphi(2^{-q}\hbox{D})v\quad\hbox{and}\quad\dot{S}_{q}u=\sum_{j\leq q-1}\dot{\Delta}_{j}u. 
$$

From the paradifferential calculus introduced by J.-M. Bony  \cite{b}   the product 
$uv$ can be formally divided into three parts as follows: 
$$
uv=T_u v+T_v u+R(u,v),
$$
where
\begin{eqnarray*}
T_u v=\sum_{q}S_{q-1}u\Delta_q v, \quad\hbox{and}\quad R(u,v)=
\sum_{q}\Delta_qu \widetilde \Delta_{q}v,
\end{eqnarray*}
$$
\textnormal{with}\quad {\widetilde \Delta}_{q}=\sum_{i=-1}^{1}\Delta_{q+i}.
$$
$T_{u}v$ is called paraproduct of $v$ by $u$ and  $R(u,v)$ the remainder term. 

 Let $(p_{1},p_{2})\in[1,+\infty]^2$ and $s\in\RR,$ then the nonhomogeneous  Besov 
\mbox{space $B_{p_{1},p_{2}}^s$} is 
the set of tempered distributions $u$ such that
$$
\|u\|_{B_{p_{1},p_{2}}^s}:=\Big( 2^{qs}
\|\Delta_q u\|_{L^{p_{1}}}\Big)_{\ell^{p_{2}}}<+\infty.
$$

Let us now recall the Lorentz spaces. For a measurable function  $f$ we define its nonincreasing rearrangement  by
$$
f^{\ast}(t) := \inf\Big\{s,\;\lambda\big(\{ x,\; |f (x)| > s\}\big)\leq t\Big\},
$$
where $\lambda$ denotes the usual Lebesgue measure.
For $(p,q)\in [0,+\infty]^2,$ the Lorentz space $L^{p,q}$ is the set of  functions $f$ 
such that $\|f\|_{L^{p,q}}<\infty,$ with
$$
\|f\|_{L^{p,q}}:=\left\lbrace
\begin{array}{l}\displaystyle
\Big(\int_0^\infty[t^{1\over p}f^{\ast}(t)]^q{dt\over t}\Big)^{1\over q}, \quad \hbox{for}\;1\leq q<\infty\\
\displaystyle\sup_{t>0}t^{1\over p}f^{\ast}(t),\quad \hbox{for}\;q=\infty.
\end{array}
\right.
$$
We can also define Lorentz spaces  by real interpolation from Lebesgue spaces:
$$
(L^{p_0},L^{p_1})_{(\theta,q)}=L^{p,q},
$$
where
$
1<p_0<p<p_1<\infty,
$
$\theta$ satisfies ${1\over p}={1-\theta\over p_0}
+{\theta\over p_1}$ and $1\leq q\leq\infty$.
We have the classical properties:
\begin{equation}\label{imbed0}
L^\infty\times L^{p,q}\rightarrow L^{p,q},
\end{equation}
\begin{equation}
\label{imbed23}L^{p,q}\hookrightarrow L^{p,q'},\forall\, 1\leq p\leq\infty; 1\leq q\leq q'\leq \infty\quad \hbox{and}\quad L^{p,p}=L^p.
\end{equation}
We will now precise the statement of Remark \ref{RK1}.
\begin{prop}\label{lorentz}
Let $1< p<3$  and $u\in B_{p,1}^{\frac{3}{p}+1}(\RR^3)$ be an axisymmetric divergence free vector field. We denote by $\omega$ its vorticity, then we have
$$
\|{\omega}/{r}\|_{L^{3,1}}\lesssim \|u\|_{B_{p,1}^{\frac{3}{p}+1}}.
$$
\end{prop}
\begin{proof}
First we intend to show that $B_{p,1}^{\frac{3}{p}-1}\hookrightarrow L^{3,1}.$ For this purpose we write
 $$
(L^{p},L^{r})_{(\theta,1)}=L^{3,1},
$$
with
$
1<p<3<r<\infty
$
and  ${1\over3}={1-\theta\over p}
+{\theta\over r}\cdot$\\
According to  Bernstein inequalities, we have
$$
B^{{3\over p}-{3\over r}}_{p,1}
\hookrightarrow B^0_{r,1}\hookrightarrow L^{r}
\quad\mbox{and}\quad
B^0_{p,1}\hookrightarrow L^{p}.
$$
Consequently
$$
(B^0_{p,1},B^{{3\over p}-{3\over r}}_{p,1})_{(\theta,1)}
\hookrightarrow L^{3,1}
$$
On the other hand we have  (see for instance \cite{BE} page 152), 
$$
(B^0_{p,1},B^{{3\over p}-{3\over r}}_{p,1})_{(\theta,1)}
=B^{\theta({3\over p}-{3\over r})}_{p,1}
=B^{{3\over p}-1}_{p,1}.
$$
This completes the proof of the imbedding result.
 From this we have
\begin{eqnarray*}
\|\nabla \omega\|_{L^{3,1}}&\lesssim& \|\nabla\omega\|_{B_{p,1}^{\frac{3}{p}-1}}
\\
&\lesssim&\|u\|_{B_{p,1}^{\frac{3}{p}+1}}.
\end{eqnarray*}
It remains now to show the following estimate
$$
\|\omega/r\|_{L^{3,1}}\lesssim\|\nabla\omega\|_{L^{3,1}}.
$$
From Proposition \ref{mim} we have  $\omega(0,0,z)=0.$ Thus we get in view of Taylor formula
$$
\omega(x_{1},x_{2},z)=\int_{0}^1\Big(x_{1}\partial_{x_{1}}\omega(\tau x_{1},\tau x_{2},z)+x_{2}\partial_{x_{2}}\omega(\tau x_{1},\tau x_{2},z)\Big)d\tau
$$
Therefore we obtain from (\ref{imbed0}) and by homogeneity
\begin{eqnarray*}
\|\omega/r\|_{L^{3,1}}&\lesssim&\int_{0}^1\|\nabla\omega(\tau\cdot,\tau\cdot,\cdot)\|_{L^{3,1}}d\tau\\
&\lesssim& \|\nabla\omega\|_{L^{3,1}}\int_{0}^1\tau^{-\frac{2}{3}}d\tau\\
&\lesssim &\|\nabla\omega\|_{L^{3,1}}.\end{eqnarray*}
This achieves the proof.
\end{proof}
Also we need the following result.
\begin{prop}\label{LZ}
Given $(p,q)\in[1,\infty]^2$ and a divergence free vector field $u$ such that $\nabla u\in L^1_{\textnormal{loc}}(\RR_{+};\,L^\infty).$ Take a smooth solution $f$ of the transport equation 
 $$
\partial_{t}f+u\cdot\nabla f=0,\, f_{|t=0}=f_0.
$$
Then we have 
$$
\|f(t)\|_{L^{p,q}}\leq\|f_{0}\|_{L^{p,q}}.
$$
\end{prop}
\begin{proof}
We use the conservation of  the $L^p$ norm combined with a standard interpolation argument.
\end{proof}
In the sequel we need the following proposition (see \cite{Ch1} page 66 for the proof).
 \begin{prop}\label{Lems2}
Given $0<s<1$ and $ \nabla u\in L^1_{\textnormal{loc}}(\RR_+;\, L^\infty).$
Let $f$ be a  solution of the transport equation 
 $$
\partial_{t}f+u\cdot\nabla f=g,\, f_{|t=0}=f_0,
$$
such that $f_0\in B_{\infty,\infty}^s(\RR^3)$ and 
$g\in{L^1_{\textnormal{loc}}}(\RR_{+};B_{\infty,\infty}^{s}).$     
\mbox{Then $\forall t\in\RR_{+},$}
$$
\|f(t)\|_{B_{\infty,\infty}^s}     
\leq 
Ce^{C\int_{0}^t\|\nabla u(\tau)\|_{L^\infty}d\tau}
\Big(\|f_0\|_{B_{\infty,\infty}^s}+\int_{0}^t\|g(\tau)\|_{B_{\infty,\infty}^{s}}d\tau\Big),
$$
where $C$ is a constant depending only on $s.$ 
\end{prop}
Let us now give some precise results concerning the  axisymmetric vector fields.
\begin{prop}\label{mim} Let $u=(u^1,u^2,u^3)$ be a smooth axisymmetric vector field. 
Then we have
\begin{itemize}
\item[(i)] for every $q\geq -1$,  $\Delta_{q}u$ is axisymmetric,
\item[(ii)] the vector $\omega=\nabla\times u=(\omega^1,\omega^2,\omega^3)$ satsifies:
$\omega^3=0$ and
$$
x_{1}\omega^1(x_{1},x_{2},z)+x_{2}\omega^2(x_{1},x_{2},z)=0,
$$
for every $(x_{1},x_{2},z)\in \Bbb R^3$,
\item[(ii)] $\omega^1(x_{1},0,z)=\omega^2(0,x_{2},z)=0$ for every $(x_{1},x_{2},z)\in \Bbb R^3$.
\end{itemize}
\end{prop}

\section{Proof of Theorem \ref{thm0}}
It is well-known from the general theory of hyperbolic systems that to prove Theorem \ref{thm0}  it is enough to give a global {\it a priori} estimate of the Lipschitz norm of the velocity. This  allows, in particular, to propagate the initial Besov regularity  $B_{p,1}^{\frac{3}{p}+1}$, see for instance \cite{Chae, Hee}. Thus we  restrict our attention to the establishment of a   Lipschitz estimate which will be done in several propositions. 
The first one gives an {\it a priori}  bound to the $L^\infty$-norm of $u^r\over r\cdot$.
\begin{prop} 
\label{u/r}
For every $t\geq 0$,
\begin{eqnarray*}
\nonumber
\big\|{u^r(t)/r}\big\|_{L^\infty}
&\lesssim&
\big\|{\omega_0/r}\big\|_{L^{3,1}}.
\end{eqnarray*}
\end{prop}
\begin{proof}
According to  Lemma 1 in \cite{Taira}, one has for every 
  $t\geq 0$
$$
|u^r(t,x)|
\lesssim \int_{| y-x|\leq r}\frac{|\omega(t,y)|}{|x-y|^{2}}dy
+r \int_{| y-x|\geq r}\frac{|\omega(t,y)|}{|x-y|^{3}}dy,
$$
with $r=\sqrt{x_{1}^2+x_{2}^2}$.
Thus, one can estimate
$$
\begin{aligned}
|u^r(t,x)|
&\lesssim
\int_{| y-x|\leq r}{\frac{|\omega(t,y)|}{r'}}{r'\over |x-y|^2}dy
+r \int_{| y-x|\geq r}{|\omega(t,y)|\over r'}{r'\over |x-y|^3}dy
\\&
\lesssim
r\int_{|y-x|\leq r}{|\omega(t,y)|\over r'}{1\over |x-y|^2}dy
+r \int_{| y-x|\geq r}{|\omega(t,y)|\over r'}{r'-r+r\over |x-y|^3}dy,
\end{aligned}
$$
where we have used the notation $r'=\sqrt{y_1^2+y_2^2}$.

  Since $\vert r'-r\vert\leq |x-y|$, then we get easily
$$
\begin{aligned}
|u^r(t,x)|
&\lesssim
r\int_{\RR^3}{|\omega(t,y)|\over r'}{1\over |x-y|^2}dy
+r^2\int_{| y-x|\geq r}{|\omega(t,y)|\over r'}{1\over |x-y|^3}dy
\\&\lesssim r\int_{\RR^3}{|\omega(t,y)|\over r'}{1\over |x-y|^2}dy
\\&
:=\hbox{I}.
\end{aligned}
$$
 As  ${1\over |\cdot|^2}\in L^{{3\over 2},\infty}(\RR^3)$, 
 then Young inequalities on  $L^{q,p}$ spaces
 \footnote{The convolution 
$L^{p,q}\star L^{p',q'} \longrightarrow L^{\infty}$
is a bilinear  continuous operator, (see \cite{Neil}, page 141 for more details).} 
  imply
$$
\hbox{I}\lesssim
r\big\|{\omega / r}\big\|_{L^{3,1}}.
$$
Thus we obtain
$$
\big\|{u^r/ r}\big\|_{L^\infty}
\lesssim \big\|{\omega/ r}\big\|_{L^{3,1}}.
$$
Since  $\omega/ r$  satisfies \eqref{equation_importante} then applying Proposition \ref{LZ} concludes the proof.
\end{proof}
%%%%%%%%%%%%%%%%%%%%%%%%%%%%%%%%%%%%
Let us now show how to use this for bounding  both vorticity and  velocity.
\begin{prop}
\label{ur0}  There exist two functions $\Phi_1$ and $\Phi_2$ such that 
$$
\|\omega(t)\|_{L^\infty}\leq \Phi_1(t)
$$
and 
$$
\|u(t)\|_{L^\infty}
\leq \Phi_2(t),
$$
for all $t\geq 0$.
\end{prop}
%%%%%%%%%%%%%%%%%%%%%%%%%%%%%%%%%%%%
\begin{proof}
 From the maximum principle applied to  \eqref{tourbillon}  one has
$$
\|\omega(t)\|_{L^\infty}
\leq
\|\omega_0\|_{L^\infty}
+\int_0^t\big\|{u^r(\tau)/ r}\big\|_{L^\infty}\|\omega(\tau)\|_{L^\infty}d\tau.
$$
Using  Gronwall's lemma and Proposition \ref{u/r} gives the first bound.

To estimate $L^\infty$ norm of the velocity we write
$$
\begin{aligned}
\|u(t)\|_{L^\infty}
&\leq \|\dot S_{-N}u\|_{L^\infty}+
\sum_{q\geq - N}\|\dot \Delta_qu\|_{L^\infty},
\end{aligned}
$$
where $N$ is an arbitrary positive integer that will be fixed later.

By Bernstein inequality we infer
 $$
\sum_{q> -N}\|\dot\Delta_q u\|_{L^\infty}
\lesssim 
2^{N}\|\omega\|_{L^\infty}.
$$
For the other part we use the intergal equation  to get 
$$
\begin{aligned}
\|\dot S_{-N} u\|_{L^\infty}
&\leq \|\dot S_{-N}u_0\|_{L^\infty}+
\int_0^t\|\dot S_{-N}\big(\mathbb{P}(u\cdot\nabla)u\big)\|_{L^\infty}d\tau
\\&
\lesssim
\| u_0\|_{L^\infty}
+2^{-N}\int_0^t\|u(\tau)\|_{L^\infty}^2d\tau,
\end{aligned}
$$
where $\mathbb{P}$ denotes the Leray's projector over divergence free vector fields. Hence we obtain
$$
\|u(t)\|_{L^\infty}\lesssim\|u_{0}\|_{L^\infty}+2^N\|\omega(t)\|_{L^\infty}+2^{-N}\int_{0}^t\|u(\tau)\|_{L^\infty}^2d\tau.
$$
If we choose $N$ such that   
$$
2^{2N}\approx 1+\frac{\displaystyle \int_{0}^t\|u(\tau)\|_{L^\infty}^2d\tau}{\|\omega(t)\|_{L^\infty}},
$$
then we obtain
$$
\|u(t)\|_{L^\infty}^2
\lesssim
\| u_0\|_{L^\infty}^2+\|\omega(t)\|_{L^\infty}^2
+\|\omega(t)\|_{L^\infty}
\int_0^t\|u(\tau)\|_{L^\infty}^2d\tau.
$$
Thus Gronwall's  lemma and the boundedness of the vorticity yield
\begin{equation}
\begin{aligned}
\label{u_infty}
\|u(t)\|_{L^\infty}
&\lesssim
\big(\|u_0\|_{L^\infty}+\|\omega\|_{L^\infty_{t}L^\infty}\big)e^{Ct\|\omega\|_{L^\infty_{t}L^\infty}}
\\&
\leq
\Phi_{2}(t).
\end{aligned}
\end{equation}
\end{proof}
In the next proposition we give   some  precise estimates on the velocity.
%%%%%%%%%%%%%%%%%%%%%%%%%%%%%%%%%
\begin{prop}
\label{ur1}  There exists a function $\Phi_3$  such that for all $t\in\RR_{+},$
$$
\|\omega(t)\|_{B^0_{\infty,1}}+\| u(t)\|_{B_{\infty,1}^1}
\leq \Phi_3(t).
$$
\end{prop}
\begin{proof}
We will use for this purpose a new approach similar to \cite{ST}. 

Let $q\geq -1$ and  denote by $\tilde{\omega}_{q}$  the unique solution of the IVP
$$
\left\lbrace
\begin{array}{l}
\partial_t \tilde{\omega}_{q}+(u\cdot\nabla)\tilde{\omega}_{q}=u^r{{\tilde{\omega}_q}\over r}\\
{\tilde\omega_q}{_{|t=0}}=\Delta_{q}\omega_{0}.
\end{array}
\right.
$$
By linearity and uniqueness on has
\begin{equation}\label{LU}
\omega(t,x)=\sum_{q\geq-1}\tilde{\omega}_{q}(t,x).
\end{equation}
Using the maximum principle and Proposition \ref{u/r} we get 
\begin{eqnarray}
\nonumber
\|\tilde \omega_{q}(t)\|_{L^\infty}
&\leq& \|\Delta_{q}\omega_{0}\|_{L^\infty}e^{\int_0^t 
\|{u^r(\tau)/r}\|_{L^\infty}d\tau}
\\
\label{ben0}
&\leq&\|\Delta_{q}\omega_{0}\|_{L^\infty}\Phi_1(t).
\end{eqnarray}
We set $\mathcal{R}_{j}(t,x):=\sum_{q\geq j}\tilde \omega_{q}(t,x).$ It is obvious that
$$
\partial_t \mathcal{R}_{j}+(u\cdot\nabla)\mathcal{R}_{j}=\frac{u^r}{r}{{\mathcal{R}_j}}$$
Again from maximum principle and Proposition \ref{u/r} we get 
\begin{eqnarray}
\nonumber
\|\mathcal{R}_{j}(t)\|_{L^\infty}
&\leq& \|\mathcal{R}_{j}(0)\|_{L^\infty}e^{\int_0^t 
\|{u^r(\tau)/r}\|_{L^\infty}d\tau}
\\
\label{ben}
&\leq&\|\mathcal{R}_{j}(0)\|_{L^\infty}\Phi_1(t).
\end{eqnarray}
Also since the $z$-component of $\omega_{0}$ is zero then $\tilde{\omega}_{q}=(\tilde{\omega}_{q}^1,\tilde{\omega}_{q}^2,0)$.  We are going to work with the two components separately. The analysis will be exactly the same  so we deal only with the frist component $\tilde{\omega}_{q}^1$.
From the identity ${{u^r}\over{ r}}={{u^1}\over{x_{1}}}={{u^2}\over{ x_{2}}},$ which is an easy consequence of $u^\theta=0$, it is plain that 
the  \mbox{functions $\tilde{\omega}_{q}^1$} is solution of 
$$
\left\lbrace
\begin{array}{l}
\partial_t \tilde{\omega}_{q}^1+(u\cdot\nabla)\tilde{\omega}_{q}^1=u^2{{\tilde{\omega}_q^1}\over x_{2}},\\
{\tilde\omega_q}^1{_{|t=0}}=\Delta_{q}\omega_{0}^1.
\end{array}
\right.
$$
Let $0<\varepsilon<1,$ then by applying Proposition \ref{Lems2} to this equation  and using the fact that the H\"older sapce 
$B_{\infty,\infty}^\varepsilon:=C^\varepsilon$ is a Banach algebra, we obtain   
\begin{eqnarray}\label{C1_*}
\nonumber\|\tilde{\omega}_{q}^1(t)\|_{C^\varepsilon}
&\lesssim&
\Big(\|\Delta_{q}{\omega}_0^1\|_{C^\varepsilon}
+\int_0^t\big\|u^2{(\tilde{\omega}_{q}^1/x_{2})}(\tau)\big\|_{C^\varepsilon}d\tau\Big)
e^{CU(t)}
\\&
\lesssim&
\Big(\|\Delta_{q}{\omega}_0^1\|_{C^\varepsilon}
+\int_0^t\|u(\tau)\|_{C^{\varepsilon}}\big\|{\tilde{\omega}_{q}^1(\tau)/ x_{2}}\big\|_{C^\varepsilon}d\tau\Big)
e^{CU(t)},
\end{eqnarray}
where $U(t):=\displaystyle \int_{0}^t\|\nabla u(\tau)\|_{L^\infty}d\tau.$

 It is easy to check that  $\tilde{\omega}_{q}^1/x_{2}$ is advected by the flow, that is
$$
\left\lbrace
\begin{array}{l}
\partial_t {\tilde{\omega}_{q}^1\over x_{2}}+(u\cdot\nabla){\tilde{\omega}_{q}^1\over x_{2}}=0\\
{{\tilde{\omega}_{q}\over x_{2}}}^1{_{|t=0}}={\Delta_{q}\omega_{0}^1\over x_{2}}.
\end{array}
\right.
$$
Thus we  deduce from the maximum principle  and Proposition \ref{Lems2}, 
$$
\big\|{\tilde{\omega}_{q}^1(t)/ x_{2}}\big\|_{L^\infty}
\leq
\big\|{\Delta_q{\omega_0}^1/x_{2}}\big\|_{L^\infty}
$$
and
\begin{equation}\label{ZF1}
\big\|{\tilde{\omega}_{q}^1(t)/x_{2}}\big\|_{C^\varepsilon}
\lesssim \big\|{\Delta_q{\omega_0}^1/x_{2}}\big\|_{C^\varepsilon}
e^{CU(t)}.
\end{equation}
%%%%%%%%%%%%%%%%%%%%%%%%%%%%%%%%%%
At this stage we need the following lemma.
\begin{lem} Under the assumptions on $u_0$, one has 
$$
\Big\|{\Delta_q{\omega_0}^1/x_{2}}\Big\|_{C^\varepsilon}
\lesssim
2^{q(1+\varepsilon)}\|\Delta_q\omega_0\|_{L^\infty},\qquad \forall\, \varepsilon\in ]0,1[ .
$$
\end{lem}
\begin{proof}
Since $u_0$ is axisymmetric then according to Proposition \ref{mim}, $\Delta_q u_0$  is  too. Consequently $\Delta_q\omega_0$ is the rotationnel of an  axisymmetric vector field and then  by Proposition \ref{mim} and Taylor expansion
$$
\Delta_q{\omega_0}^1(x_{1},x_{2},z)=x_{2}\int^1_{0}(\partial_{x_{2}} \Delta_q{\omega_0}^1)(x_{1},\tau x_{2},z)d\tau.
$$
On the other hand, we have the classical equivalence norm
$$
\|v\|_{C^\varepsilon}\approx\|v\|_{L^\infty}+\sup_{x\neq y}\frac{|v(x)-v(y)|}{|x-y|^\varepsilon}\cdot
$$
This leads to
$$
\|v(\cdot,\tau\cdot,\cdot)\|_{C^\varepsilon}\lesssim \|v\|_{C^\epsilon},\quad\forall\tau\in [0,1].
$$
 Hence we get
 \begin{eqnarray*}
\big\|{\Delta_q{\omega_0}^1/x_{2}}\big\|_{C^\varepsilon}
&\leq& \int^1_{0} \Vert( \partial_{x_{2}}\Delta_q{\omega_0}^1)(\cdot,\tau\cdot,\cdot)\Vert_{C^\varepsilon}d\tau
\\
&\leq&\|\partial_{x_{2}}\Delta_q\omega_0^1\|_{C^\varepsilon}
\\
&\lesssim& 
2^{q(1+\varepsilon)}\| \Delta_q\omega_0^1\|_{L^\infty},
\end{eqnarray*}
as claimed.
\end{proof}
From (\ref{ZF1}) and the above lemma  we infer
\begin{equation}
\nonumber
\big\|{\tilde{\omega}_{q}^1/x_{2}}(t)\big\|_{C^\varepsilon}
\lesssim2^{q(1+\varepsilon)}\|\Delta_q\omega_0^1\|_{L^\infty}
e^{CU(t)}.
\end{equation}
Plugging this estimate into (\ref{C1_*})
\begin{equation}\label{ZF2}
\|\tilde{\omega}_{q}^1(t)\|_{C^\varepsilon}
\lesssim
2^{q(1+\varepsilon)}\|\Delta_{q}{\omega}_0^1\|_{L^\infty}
\Big(1+\int_0^t\|u(\tau)\|_{C^\varepsilon}d\tau\Big)
e^{CU(t)}
\end{equation}
It remains to estimate $\|u\|_{C^\varepsilon}.$ For this purpose we recall the classical estimates
\begin{equation*}
\|u\|_{C^\varepsilon}\lesssim \|u\|_{B^1_{\infty,\infty}}
\lesssim
\|u\|_{L^\infty}+\|\omega\|_{L^\infty}.
\end{equation*}
Now according to Proposition \ref{ur0}, we have
\begin{equation}
\label{u_Besov}
\|u\|_{C^\varepsilon}\leq \Phi_2(t).
\end{equation}
Combining  (\ref{ZF2}) and  (\ref{u_Besov}) with de definition of Besov spaces we get  for all $j\geq-1$
\begin{eqnarray}
\nonumber
\|\Delta_{j}\tilde{\omega}_{q}^1(t)\|_{L^\infty}
&\lesssim&  2^{-\varepsilon j+{(1+\varepsilon)q}}\|\Delta_q\omega_0^1\|_{L^\infty}\Phi_2(t)e^{CU(t)}
\\\nonumber
&\lesssim& 2^{-\varepsilon(j-{1+\varepsilon\over\varepsilon}q)}\|\Delta_q\omega_0\|_{L^\infty}\Phi_2(t)e^{CU(t)}.
\end{eqnarray}
Similar arguments give the same estimate for $\tilde{\omega}_{q}^2$. Thus,  we obtain  \mbox{ for $\varepsilon=\frac{1}{2}$}
\begin{equation}\label{t1}
\|\Delta_{j}\tilde{\omega}_{q}(t)\|_{L^\infty}  \lesssim 2^{-\frac{1}{2}(j-3q)}\|\Delta_q\omega_0\|_{L^\infty}\Phi_2(t)e^{CU(t)},\quad\forall\, j,q\geq-1.
\end{equation}
Let $N$ be a fixed positive integer that will be carefully chosen later. Then   we have from (\ref{LU})
\begin{eqnarray}
\nonumber
\|\omega(t)\|_{B_{\infty,1}^0}
&\leq&
\sum_{j}\|\Delta_{j}\sum_{q}\tilde{\omega}_{q}(t)\|_{L^\infty}
\\\nonumber
&\leq&
\sum_{j-3q\geq N}\|\Delta_{j}\tilde{\omega}_{q}(t)\|_{L^\infty}+\sum_{3\leq j-3q<N}\|\Delta_{j}\tilde{\omega}_{q}(t)\|_{L^\infty}\\
\nonumber&+&\sum_{j}\|\Delta_{j}\sum_{j-3<3q}\tilde{\omega}_{q}(t)\|_{L^\infty}
\\
\label{t2}
&:=&\hbox{I}+\hbox{II}+\hbox{III}.
\end{eqnarray}

To estimate the first term we use (\ref{t1})  and  the convolution inequality for the series
\begin{equation}
\label{t3}
\hbox{I}
\lesssim
2^{-{1\over 2} N}\|\omega^0\|_{B_{\infty,1}^0}e^{CU(t)}\Phi_2(t).
\end{equation}

To estimate the term $\hbox{II}$ we use two facts: the first
 one is that the operator $\Delta_{j}$ maps uniformly  $L^\infty$ into itself while the second  is  the $L^\infty$ 
 estimate (\ref{ben0}),
\begin{equation}\label{ben55}
\begin{aligned}
\hbox{II}
&\lesssim
\sum_{3\leq j-3q< N}\|\tilde{\omega}_{q}(t)\|_{L^\infty}
\\&
\lesssim\Phi_1(t)\sum_{3\leq j-3q< N}\|\Delta_q\omega_0\|_{L^\infty}
\\&
 \lesssim \Phi_1(t) N\sum_{q}\|\Delta_q\omega_0\|_{L^\infty}
\\&
\lesssim N\Phi_1(t).
\end{aligned}
\end{equation}
For the third term we write in view of (\ref{ben})
\begin{eqnarray*}
\hbox{III}&\lesssim&\sum_{j}\|\sum_{q> \frac{j}{3}-1}\tilde\omega_{q}(t)\|_{L^\infty}\\
&\lesssim&\sum_{j}\|\mathcal{R}_{[\frac{j}{3}]}(t)\|_{L^\infty}\\
&\lesssim&\Phi_{1}(t)\sum_{j}\|\mathcal{R}_{[\frac{j}{3}]}(0)\|_{L^\infty}\\
&\lesssim&\Phi_{1}(t)\|\omega_{0}\|_{\widetilde B_{\infty,1}^0}.
\end{eqnarray*}
Combining this estimate with (\ref{ben55}),  (\ref{t3}) and (\ref{t2}) we obtain
$$
\begin{aligned}
\|\omega(t)\|_{B_{\infty,1}^0}
&\lesssim 2^{-{1\over 2} N}e^{CU(t)}\Phi_2(t)+  N\Phi_1(t)
\\&
\lesssim (2^{-{1\over 2} N}e^{CU(t)}+N)\Phi_2(t).
\end{aligned}
$$
Putting
$$
N=\Big[\frac{2CU(t)}{\log 2}\Big]+1,
$$
we obtain
$$
\|\omega(t)\|_{B_{\infty,1}^0}
\lesssim \big(U(t)+1\big)\Phi_2(t).
$$
On the other hand we have
$$
\|\nabla u\|_{L^\infty}
\lesssim
\|u\|_{L^\infty}+\|\omega\|_{B^0_{\infty,1}},
$$
which yields, via Proposition \ref{ur0},
$$
\begin{aligned}
U(t)&\lesssim \int^t_0\|u(\tau)\|_{L^\infty}d\tau+ \int^t_0\|\omega(\tau)\|_{B^0_{\infty,1}}d\tau
\\&
\lesssim \Phi_2(t)+ \int^t_0\|\omega(\tau)\|_{B^0_{\infty,1}}d\tau.
\end{aligned}
$$
Hence we obtain
$$
\|\omega(t)\|_{B_{\infty,1}^0}
\leq \Phi_2(t)\big( 1+ \int^t_0\|\omega(\tau)\|_{B^0_{\infty,1}}d\tau\big).
$$
It follows from Gronwall's lemma that
$$
\|\omega(t)\|_{B_{\infty,1}^0}
\leq \Phi_3(t).
$$
We have also obtained the estimate
$$
\|\nabla u\|_{L^\infty}+\|u\|_{B_{\infty,1}^1}
\leq\Phi_3(t).
$$
This concludes the proof of Proposition \ref{ur1}.
\end{proof}

\end{document}